\documentclass[12pt]{amsart}

\input{xy}
\xyoption{all}
\usepackage{graphicx}
\usepackage{color}
\usepackage{amssymb}
\usepackage{amsmath}
\usepackage{graphics}

\newtheorem{theorem}{Theorem} [section]
\newtheorem{prop}[theorem]{Proposition}
\newtheorem{lemma}[theorem]{Lemma}
\newtheorem{cor}[theorem]{Corollary}

\numberwithin{equation}{section}
\numberwithin{figure}{section}

\newcommand\C{{\mathbb C}}

\renewcommand\P{{\mathbb P}}
\newcommand\R{{\mathbb R}}
\newcommand\Z{{\mathbb Z}}

\newcommand\N{{\mathbb N}}
\newcommand\D{{\mathbb D}}

\renewcommand\phi{\varphi}

\newcommand\del{\partial}

\newcommand\M {\mathrm{M}}

\newcommand\Precrit {\operatorname{Precrit}}
\newcommand\Crit {\operatorname{Crit}}

\begin{document}

\title{THE SCHWARZIAN DERIVATIVE AND POLYNOMIAL ITERATION}

\author{Hexi Ye}

\begin{abstract}
We consider the Schwarzian derivative $S_f$ of a complex polynomial $f$ and
its iterates. We show that the sequence $S_{f^n}/d^{2n}$ converges to $-2(\partial G_f)^2$, for $G_f$ the escape-rate function of $f$.  As a quadratic differential, the Schwarzian derivative $S_{f^n}$ determines a conformal metric on the
plane.  We study the ultralimit of these metric spaces.

\end{abstract}

\date{\today}

\maketitle

% For editing purposes only
%\tableofcontents

% suppress page number 1
\thispagestyle{empty}

%%%%%%
%%%%%%

\section{Introduction}
Recall that the Schwarzian derivative of a holomorphic function $ f $ on the complex plane is defined as
   $$S_f(z)\ = \frac{f'''}{f'}- \frac{3}{2}\left(\frac{f''}{f'}\right)^2$$
It is well known that $S_f\equiv0$ if and only if $f$ is a M$\ddot{\textup{o}}$bius transformation. We can view the Schwarzian derivative as a measure of the complexity of a nonconstant holomorphic function.

%Our goal is to detect dynamical $f$-invariant objects that arise as the `` limit "  of $\{S_{f^n}\}$. After rescaling, we may expect that it converges (in some sense) to some thing coming from the dynamics of $f$.

%In this article, we consider the possible limit of the rescaled Schwarzian Derivative for polynomials with degree $d \geq 2$, i.e. the possible limit of $\{\frac{S_{f^n}}{d^{2n}}\}$ as $n$ tends to infinity. Hereafter, $f^n$ means $f$ composed with itself $n$ times.

%For any critical point $z_o$ of a holomorphic function f, the function $S_f(z)$ has a pole of order 2 at $z_o$. Moreover, the coefficient of  $\frac {1}{(z-z_o)^2}$ term at this pole is uniquely determined by the local degree of $f$ at this critical point $z_o$.

Let $f:\C \to \C$ be a complex polynomial with degree $d\geq 2$. In this article we examine the sequence of Schwarzian derivatives of the iterates $f^n$ ($f$ composed with itself $n$ times) of $f$. Specifically, we look at the sequence
$$ \left\{ \frac{S_{f^n}(z)}{d^{2n}}\right\}_{n\geq 1}$$
 and view it as a sequence of meromorphic functions or quadratic differentials on the Riemann sphere. We are interested in understanding the limit as $n\to \infty$.

We begin with the simplest example.

\bigskip

{\bf Example 1.}
Let $f(z)=z^d$ with $d\geq 2$, then we get
$$
\begin{array}{lll}
S_{f^n}(z)
&=\frac{2(d^n-1)(d^n-2)-3(d^n-1)^2}{2z^2}\\[6pt]
&=\frac{1-d^{2n}}{2z^2}
\end{array}
$$
Since $d\geq2$, the sequence of normalized Schwarzians converges,
$$\lim_{n\to\infty} \frac{S_{f^n}}{d^{2n}}=\lim_{n\to\infty} \frac{1-d^{2n}}{2d^{2n}z^2}= -\frac{1}{2z^2}$$
locally uniformly on $\C \backslash \{0\}$.

%We have similar limit for all the polynomials conformal conjugated to $z^d$ with $d\geq 2$.

The normalized Schwarzians as quadratic differentials $\left\{\frac{S_{f^n}}{d^{2n}}dz^2\right\}$ converge to $-\frac{1}{2z^2}dz^2$. The associated conformal metric $ds=\frac{|dz|}{|\sqrt{2}z|}$ makes $\C\backslash \{0\}$ isometric to an infinite cylinder of radius $\frac{1}{\sqrt{2}}$. The cylinder's closed geodesics are the horizontal trajectories of the quadratic differential $-\frac{1}{2z^2}dz^2$.

\bigskip

{\bf Local convergence.} Let $G_f$ be the escape-rate function of $f$, which is defined as
 $$G_f=\lim_{n\to \infty} \frac{\log^+|f^n|}{d^{n}},$$
where $\log^+|x|=\textup{max}(\log |x|,0)$. Let $\Precrit(f)=\cup_{n>0}\{c\in \C|(f^n)'(c)=0\}$ be the union of the critical points of $f$ and their backward orbits. Note that its closure $\overline{\textup{Precrit($f$)}}$ contains the Julia set $J(f)$ when $f$ is not conjugate to $z^d$.

%In section \ref{Local convergence}, we are going to prove the following convergence theorem,

\begin{theorem}\label{maintheorem1}
Let $f$ be a polynomial with degree $d\geq 2$ and not conformally conjugate to $z^d$. Then the sequence of Schwarzian derivatives $S_{f^n}$ satisfies
$$\lim_{n\to\infty}\frac{S_{f^n}(z)}{d^{2n}}=-2\left(\frac{\del G_f(z)}{\del z}\right)^2,$$
locally uniformly on $\C \backslash \overline{\Precrit(f)}$.
\end{theorem}

{\bf Remark.} The choice of normalization $\frac{1}{d^{2n}}$ allows us to focus on the basin of infinity. Other normalizations might detect interesting  properties of $J(f)$. In Corollary \ref{maincor}, we show that $\left\{\frac{S_{f^n}}{d^{2n}}dz^2\right\}$ converge on the entire Fatou set, in the sense of $L^{\frac{1}{2}}_{loc}$ convergence.

Sometimes, people are also interested in the nonlinearity of a nonconstant holomorphic function on $\C$.  Similar with Theorem \ref{maintheorem1}, we have
$\lim_{n\to \infty} \frac{N_{f^n}dz}{d^n}=\partial G_f$ locally uniformly on $\C\setminus \overline{\Precrit (f)}$, where $N_f=f'/f''$; see Theorem \ref{nonlinearity}.

\bigskip

%In section \ref{Ultralimit of Metric Spaces}.
{\bf Metric space convergence.} Let $f$ be a polynomial with degree $d\geq 2$. Each $S_{f^n}$ determines a conformal geodesic metric $d_n$ on the complement of the critical points of $f^n$, given by $ds=${\tiny $\sqrt{|-\frac{S_{f^n}}{d^{2n}}dz^2|}$}. From this sequence of geodesic spaces, we obtain an ultralimit $(X_\omega,d_\omega, a_\omega)$; see Chapter I \S 5 \cite{MA} for more details about the ultralimit. The limit space is a complete geodesic space.

The escape-rate function $G_f$ also determines a conformal metric on the basin of infinity $X_o=\{z\in \C|f^n(z)\to \infty\}$ of $f$,  given by $ds=\sqrt{2}|\partial G_f|$. Given a choice of
base point $a\in X_o\backslash \Precrit (f)$, we denote this pointed metric space by $(X_o, d_o, a)$; compare \cite{DM:trees} or \cite{DP} where this metric already appeared.
\begin{figure}
  \includegraphics[width=145mm]{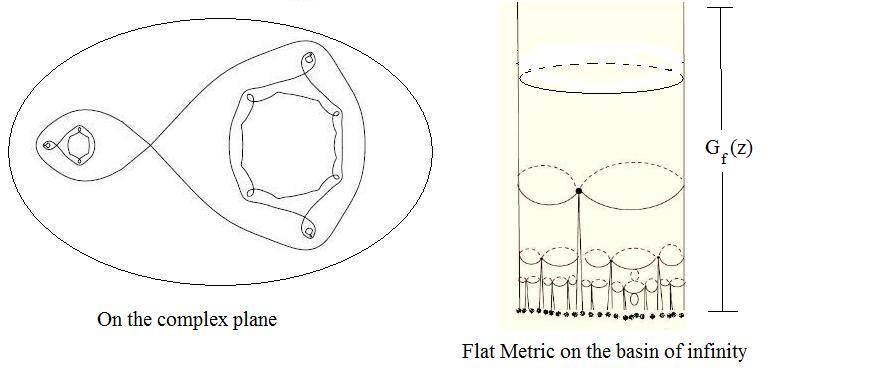}\\
    \caption{Level set structure of $G_f$ for a cubic polynomial and flat metric structure on the basin of infinity}\label{flat metric}
\end{figure}

When the Julia set is not connected, the ultralimit $X_\omega$ can be described as a ``hairy" version of $X_o$.
\begin{theorem}\label{ultralimit spaces}
 Let $f$ be a polynomial with degree $d\geq 2$ and disconnected Julia set. There exists a natural embedding from $X_o(f)\backslash \Precrit(f)$ to $X_\omega$ which extends to the metric completion $(\overline{X}_o,d_o,a)$ as an isometric embedding.

 The metric space $X_\omega$ is obtained by attaching a real ray to $\overline{X}_o$ at each point in $\Precrit(f)$, and attaching some non trivial space (containing infinitely many real rays) to each connected component of $\overline{X}_o\backslash X_o$.
\end{theorem}
{\bf Remark.} For $J(f)$ disconnected, we use the tree structure on the basin of infinity to show the embedding is an isometry. For $J(f)$ connected, we also have a locally isometric embedding of $X_o$ to $X_\omega$; this follows easily from the argument in the proof of the above theorem. But we do not expect the embedding to be a global isometry. The reason we use the ultralimit to study the limiting space is that
the spaces are not compact and the metrics $d_n$ are not uniformly
proper, so more classical notions of convergence like Gromov-Hausdorff
convergence won't work. We are not quite sure whether the ultralimit $X_\omega$ depends on the ultrafilter $\omega$ or not. But from the above theorem, the only things that might depend on the ultrafilter are the spaces attached to $\overline{X}_o\backslash X_o$.

%We associate a probability measure $\lambda_n$ to each quadratic differential $S_{f^n}dz^2$, which is determined by the length of closed trajectories near critical points of $f^n$ on $\C$.
%Let $\mu^*$ be the unique equilibrium measure of filled Julia set $K(f)$; it is also the measure of maximal entropy of $f$, see \cite{L} and \cite{FLM}. For the precise definition of $\lambda_n$, please refer to \S \ref{unitdisti}.

%\begin{theorem}\label{distribution}
%Let $f$ be a polynomial with degree $d\geq 2$ and not conformally conjugate to $z^d$. The sequence of probability measures $\{\lambda_n \}$ from the Schwarzian derivative of $f^n$ converges to $\mu^*$ weakly.
%\end{theorem}
\medskip
{\bf Conjugacy classes.} The study of $\{S_{f^n}\}$ has grown out of an attempt to better understand the moduli space $\M_d$ of polynomials (the space of conformal
conjugacy classes). The geometric structure $(f, X_o, |\partial G_f|)$ has been studied in \cite{DP} and used to classify topological
conjugacy classes. We can also use the Schwarzian derivative to classify polynomials with the same degree.  We define an equivalence relation on the set of polynomials with degree $d\geq 2$ as: $f\sim g$ if $S_fdz^2=A^*(S_gdz^2)$, for $A(z)=az+b$ some affine transformation. From this definition, polynomials $f$ and $g$ are equivalent to each other if and only if $f=B\circ g \circ A$ with $A$ and $B$ affine transformations, and if and only if $f$ and $g$ have the same critical set (counted with multiplicities) up to some affine transformation, i.e. $A(\Crit(f))=\Crit(g)$ for some affine transformation $A$. See Lemma \ref{the same critical points} for details. Note that affine conjugate polynomials are equivalent in this sense.

\begin{theorem}\label{schw equivalent}
Let $f$ and $g$ be polynomials with the same degree $d\geq 2$. Then the following are equivalent:
\begin{itemize}
\item $f^n\sim g^n$ for infinitely many $n\in \N^*$.
\item $f^n\sim g^n$ for all $n\in \N^*$.
\item $f$ and $g$ have the same Julia set up to some affine transformation ($A(J(g))=J(f)$ with A affine transformation).
\end{itemize}
Polynomials $f$ and $g$ which satisfy the above conditions are called strongly equivalent. Each such strong equivalence class consists of finitely many affine conjugacy classes (no more than the order of the symmetry group of the Julia set).
\end{theorem}

Notes: the above theorem relies on the classification of polynomials with the same Julia set, and the proof uses the main result of \cite{Bear} (see \S \ref{definition of quadratic differential} for other references and historical context).

\bigskip

I would like to thank David Dumas, Kevin Pilgrim, Stefan Wenger, Curt McMullen and especially Laura DeMarco for lots of helpful comments and suggestions.

%%%%%%
%%%%%%

\bigskip\section{Basic properties of Schwarzian derivatives }  \label{schwarzian derivative as quadratic differentials}
In this section, we give useful formulas for $S_{f^n}$ and also basic definitions that we are going to use later in this article.
 \subsection{Basic formula for the Schwarzian derivative of $f^n$}
 In order to find the limit of $\frac{S_{f^n}}{d^{2n}}$ for a general polynomial with degree $d\geq 2$, we need to rewrite $S_{f^n}$ in terms of $S_f$ and then evaluate the limit. To do this, we use the formula for  the Schwarzian derivative of the composition of two functions $f$ and $g$. An easy calculation shows that (cocycle property)
 \begin{equation}\label{basic schw dicompositon}
	S_{f \circ g}(z)=S_f(g(z))(g'(z))^2+S_g(z)
\end{equation}

From this relation we derive the following important formula:
\begin{eqnarray}\notag
S_{f^n}(z)&=&S_f(f^{n-1}(z))((f^{n-1}(z))')^2+S_{f^{n-1}}(z)\\[6pt] \label{schw equation}
&=&S_f(f^{n-1}(z))((f^{n-1}(z))')^2+S_f(f^{n-2}(z))((f^{n-2}(z))')^2+S_{f^{n-2}}(z)\\[6pt]\notag
&=&\sum_{i=1}^{n-1}S_f\circ f^{i}(z)((f^{i}(z))')^2+S_f(z)
\end{eqnarray}
\begin{prop}\label{basic prop}
Let $f$ be a polynomial with degree $d\geq 2$. For any point $z\in\C\backslash \textup{Precrit($f$)} $ and any sequence $\{n_i\}_{i=1}^\infty$ of $\N$ with $n_i\to\infty$ as $i\to\infty$, the sequences $\{\frac{S_{f^{n_i}}(z)}{d^{2n_i}}\}$ and $\{\frac{S_{f^{n_i-1}}(f(z))}{d^{2(n_i-1)}}\}$ either both converge, diverge to infinity or diverge.

Moreover, if both of them converge, then
$$d^2\lim_{i\to\infty}\frac{S_{f^{n_i}}(z)}{d^{2n_i}}=(f'(z))^2\lim_{i\to\infty}\frac{S_{f^{n_i-1}}(f(z))}{d^{2(n_i-1)}}.$$

\end{prop}

\proof From (\ref{basic schw dicompositon}), we have:
\begin{eqnarray*}
\frac{S_{f^{n_i}}(z)}{d^{2n_i}}&=&\frac{S_{f^{n_i-1}}(f(z))(f'(z))^2+S_f(z)}{d^{2n_i}}\\[6pt]
&=&\frac{S_{f^{n_i-1}}(f(z))(f'(z))^2}{d^2  d^{2(n_i-1)}}+\frac{S_f(z)}{d^{2n_i}}
\end{eqnarray*}
By the assumption that $z$ is not a critical point of $f$ and $d\geq2$, it is easy to see that this proposition is satisfied since $f'(z)$ is not equal to zero and $S_f(z)$ is finite.\qed

\medskip

\subsection{The Schwarzian derivative as a quadratic differential}\label{definition of quadratic differential}

In this subsection, we are not only considering polynomials, but also rational maps with degree $d\geq 2$.

{\bf Meromorphic  quadratic differentials.}
For any Riemann surface $S$, a meromorphic quadratic differential $Q$ on $S$ is a section of the second tensor power of the cotangent bundle. In local coordinate,
$Q=Q_1(z)dz^2$, where $Q_1(z)$ is a meromorphic function. And under changing of coordinate $w=w(z)$,
$$Q=Q_2(w)(dw)^2=Q_2(w(z))(w'(z))^2dz^2$$
i.e., $Q_2(w(z))(w'(z))^2=Q_1(z)$.

Consider a non constant holomorphic map $f: {\P}^{1} \to {\P}^{1}$. This is a rational map with finite degree. Let's look at the Schwarzian derivative of
this rational  map $f$, and view it as quadratic differential, i.e.
$$S_fdz^2\ \  \text{instead of} \ S_f(z)$$
From the definition of the Schwarzian derivative, it is not hard to show $S_g\equiv 0$ if and only if $g$ is a M$\ddot{\textup{o}}$bius transformation; see \cite{D}.
From the following identity
$$S_{f \circ g}(z)dz^2=S_f(g(z))(dg)^2+S_g(z)dz^2,$$
for any two M$\ddot{\textup{o}}$bius transformations $g_\circ, g_1$, we have
$$S_{g_1\circ f \circ g_\circ}dz^2=S_f\circ g_\circ (dg_\circ)^2$$
So the Schwarzian derivative $S_fdz^2$ as a quadratic differential is well defined on $\P^1$. More generally, for any non constant holomorphic map from some projective Riemann surface to another projective Riemann surface, there is an unique quadratic differential (named as Schwarzian derivative) associated to it; see \cite{D}.

Recall that in the last part of the introduction, we defined an equivalence relation of polynomials: $f\sim g$ if the Schwarzian derivative of $f$ is the same as the Schwarzian derivative of $g$ up to some affine transformation. The Schwarzian derivative of a polynomial is determined by the locations and multiplicities of the critical points:
\begin{lemma}\label{the same critical points}
Let $f$ and $g$ be polynomials with degree $d\geq 2$.  Then $f\sim g$ if and only if they have the same critical set (critical points are counted with multiplicity) up to some affine transformation ($A(\Crit(f))=\Crit(g)$ with $A$ some affine transformation).
\end{lemma}
\proof Assume $f\sim g$, then $S_fdz^2=A^*(S_gdz^2)$, which means $f=B\circ g \circ A$ for some affine transformations $A$ and $B$. Indeed, by the cocycle property,  $A^*(S_gdz^2)=S_{g\circ A}dz^2$ and then $S_{f\circ(g\circ A)^{-1}}\equiv0$ on some open subset of $\C$. So $B=f\circ(g\circ A)^{-1}$ is a M$\ddot{\textup{o}}$bius transformation on some open subset of $\C$. By continuity, $f=B\circ g \circ A$ in $\C$, and so $B$ is an affine transformation. This  implies that $A$ transforms the critical set of $f$ to the critical set of $g$.

Conversely, assume that there is an affine transformation $A$ that transforms the critical set of $f$ to the critical set of $g$. Since $g\sim g\circ A$, it suffices to show that $f\sim g\circ A$. Because $f$ and $g\circ A$ have the same critical set, so we can let the critical set be $\{c_i\}_{i=1}^{d-1}$. Then $f=ah(z)+b$ and $g\circ A=ch(z)+d$ with $a, c\neq 0$ and $h(z)=\int_0^z\prod_{i=1}^{d-1}(t-c_i)dt$. Which means $S_f=S_h=S_{g\circ A}$, i.e. $f\sim g\circ A \sim g$.\qed

\bigskip
{\bf Proof of Theorem \ref{schw equivalent}}. Let $\simeq$ be the notion of strong equivalence. Assume that the polynomial $f$ with degree $d\geq 2$ is not conjugate to $z^d$, and there is a subsequence $\{n_i\}_{i=1}^\infty \subset \N^*$ such that $f^{n_i}\sim g^{n_i}$. By Lemma \ref{the same critical points}, there are affine transformations $\{A_i=a_iz+b_i\}$ such that $A_i(\Crit(f^{n_i}))=\Crit(g^{n_i})$. Since $f$ is not conjugate to $z^d$, so $g$ is not conjugate to $z^d$. Indeed, for $n\geq 2$, $\Crit(f^{n})$ has at least two distinct points, however, $\Crit(z^{d^n})$  has only one point. Set $M_1=\textup{Diam}(\Crit(f^2))>0$, $M_1'=\textup{Diam}(\Crit(g^2))>0$, $M_2=\textup{Diam}(\Precrit(f))$ and $M_2'=\textup{Diam}(\Precrit(g))$. Because $f(\C\backslash D(0,R))\subset \C\backslash D(0,R)$ for $R$ sufficiently large, so $\Precrit(f)$ is bounded and then $M_2<\infty$. Similarly,  $M_2'<\infty$. Moreover, since $f^{n+m}=f^n\circ f^m$, then $\Crit(f^{m})\subset \Crit(f^{n+m})$. So for the diameters Diam($\Crit(f^{n_i})$) and Diam($\Crit(g^n_i)$) of the critical sets, we have
$$ 0< M_1\leq Diam(\Crit(f^{n_i}))\leq M_2<\infty,\ 0< M_1'\leq Diam(\Crit(g^{n_i}))\leq M_2'<\infty,$$
for any $n_i\geq 2$. Consequently,
$$0<M_3\leq |a_i| \leq M_4\leq \infty, |b_i|\leq M_5<\infty$$
So after passing to a subsequence and without loss of generality, we can assume $A_i\to A$ as $i\to \infty$, where $A$ is an affine transformation. Then it is easy to know $A(J(f))=J(g)$. Indeed, for any $c\in \Crit(f^j)\subset \Crit(f^i) $, $A_i(c)\in \Crit(g^i)\subset \Precrit(g)$ with $j\leq i$. It indicates that $A(c)\in \overline{\Precrit(g)}$ and then $A(\overline{\Precrit(f)})\subset  \overline{\Precrit(g)}$. Similarly, by taking $A_i^{-1}$ instead of $A_i$, we get $A^{-1}(\overline{\Precrit(g)})\subset \overline{\Precrit(f)}$. Because $f$ is not conjugate to $z^d$, the set of accumulating points of $\overline{\Precrit(f)}$ (respt. $\overline{\Precrit(g)}$) is $J(f)$ (respt. $J(g)$), which means that $A(J(f))=J(g)$.

If $f$ is conjugate to $z^d$, then by the above argument, $g$ should also be conjugate to $z^d$. So there is an affine map $A$ such that $A(J(f))=J(g)$.

Conversely, assume $f$ and $g$ have the same degree $d\geq 2$ and $A(J(f))=J(g)$ for some affine transformation $A$. Since $g_1=A^{-1}\circ g\circ A\simeq g$ and $A(J(g_1))=J(g)=A(J(f))$, so it is enough to prove that $g_1 \simeq f$ with the condition that they have the same Julia set. First, if the Julia set is a circle, then both of them are conjugate to $z^d$. Indeed, we can assume $J(f)$ is the closed unit disk. Let $\phi$ be the Boettcher function of $f$, such that $\phi \circ f\circ \phi^{-1}=z^d$ on the basin of infinity. Since the basin of infinity is the complement of the unit disk and $\phi$ is a conformal map that fixs the infinity, so $\phi$ should be a rotation. Then $f$ is conjugate to $z^d$. So $f^n\sim z^{d^n}\sim g_1^n$ for any $n\in \N^*$. Second, since both $f^n$ and $g_1^n$ have the same degree and the same Julia set which is not a circle, then $f^n=\sigma_n \circ g_1^n$ with $\sigma_n$ an affine transformation in the symmetry group of the Julia set; see \cite{Bear} for details. So $S_{f^n}=S_{g_1^n}$ for any $n\in \N^*$. And moreover, when Julia set is not a circle, the order of symmetry group of the Julia set is finite; see Lemma 4 in \cite{BA}. Thus there are only finitely many conjugacy classes which belongs to a strong equivalence class.\qed

Notes: the proof of Theorem \ref{schw equivalent} relies on the classification of polynomials with the same Julia set. When do two polynomials have the same Julia set? Historically, commuting polynomials have the same Julia set, as observed by by Julia in 1922 \cite{Julia}. Later in 1987, Baker and Eremenko showed when the symmetric group of the Julia set is trivial, polynomials with this Julia set commute; see \cite{Baker}. In 1989, Fern$\acute{\textup{a}}$ndez showed there is at most one polynomial with given degree, leading coefficient and Julia set; see \cite{Fern}. Finally in 1992, Beardon showed $ \{ g | \textup{deg}(f) = \textup{deg}(g), J(f) = J(g) \} = \{\sigma \circ f | \sigma \in \textup{symmetric group of }J(f)\}$; see \cite{Bear}.

\medskip
\subsection
{\bf Conformal metric of quadratic differential.} For meromorphic quadratic differential $Q$ on $S$, it determines a flat metric $ds^2=|Q|$, with singularities at zeros and poles of $Q$.

\medskip
{\bf Trajectories as a foliation.} For any meromorphic quadratic differential $Q$ on $S$, it determines a foliation structure on $S$ with singularities at zeros and poles of $Q$. A smooth curve on $S$ is a (horizontal) trajectory of $Q$, if it does not pass though any zero or pole of $Q$, and for any point $p$ in the curve, the non zero vector $dz$ tangent to this curve satisfy:
$$\arg(Q(p)dz^2)=0$$
i.e $Q(p)dz^2$ is a positive real number. By a trajectory, we usually mean the trajectory that it is not properly contained in another trajectory, i.e.  a maximal trajectory.

\begin{lemma}\label{length of normalized trajectory}
Let $f: {\P}^{1} \to {\P}^{1}$ be a rational map with degree $d\geq 2$. Then $S_fdz^2$ is a meromorphic quadratic differential with poles of order two at the critical points of $f$. For any critical point $p$ of $f$ with order $k$, near $p$
$$S_f(z)=\frac{1-k^2}{2(z-p)^2}+O(\frac{1}{|z-p|})$$
i.e., a neighborhood of $p$ is an infinite cylinder with closed geodesics as trajectories of length $2\pi \sqrt{\frac{k^2-1}{2}}$.
\end{lemma}
\proof The only thing we need to show here is that the coefficient of $\frac {1}{(z-p)^2}$ at $p$ is $\frac {1-k^2}{2}$, and for other details of this lemma, please refer to the \S 6.3 \cite{KS}. The coefficient can be verified by a direct computation.\qed

\bigskip
%From Lemma \ref{length of normalized trajectory}, one may ask whether the quadratic differential $S_{f}dz^2$ is a Jenkins-Strebel differential, i.e., the horizontal trajectories are closed except the trajectories terminating at the singular points? Well, we are not sure whether this statement is true or not, but it is not hard to see that for all the polynomials with the number of distinct critical points less than three, and for all the polynomials with real coefficients and three distinct critical points, their Schwarzian derivatives are Jenkins-Strebel differentials.

\subsection{$L_{loc}^{\frac{1}{2}}$ integrability of quadratic differential.} For any Riemann surface $S$, we consider the space $MQ(S)$ of meromorphic quadratic differentials on $S$ .  For any $\alpha \in MQ(S)$, we say that it is $L_{loc}^{\frac{1}{2}}$ integrable, if for any point $q\in S$, there is a local coordinate at $q$, and write $\alpha$ as $h(z)dz^2$ in this coordinate, such that the integration $\int\int\sqrt{|h(z)|}dxdy$ over some neighborhood of $p$  is finite. We say that $\{\alpha_n\}\subset MQ(S)$ $L_{loc}^{\frac{1}{2}}$-converge to $\alpha \in MQ(S)$, if both $\alpha$ and $\alpha_n$ are $L_{loc}^{\frac{1}{2}}$ integrable, and for any point $p$, there is some local coordinate at $q$, and write $\alpha$ and $\alpha_n$ as $h(z)dz^2$ and $h_n(z)dz^2$ in this local coordinate, such that the integral $\int\int\sqrt{|h_n(z)-h(z)|}dxdy$ over some neighborhood of $p$ converges to $0$. Actually, the $L_{loc}^{\frac{1}{2}}$ integrable subset of $MQ(S)$ is a vector space.

Any meromorphic quadratic differential $\alpha \in MQ(S)$ with poles of order at most two is $L_{loc}^{\frac{1}{2}}$ integrable.

\begin{lemma}\label{L-1/2 integrable}
For any rational function $f$ with degree $d\geq 2$, $S_{f^n}dz^2$ is $L_{loc}^{\frac{1}{2}}$ integrable.
\end{lemma}

\proof Since for any rational map $f$, the Schwarzian derivative $S_fdz^2$ has poles of order at most two, then $S_fdz^2$ is $L_{loc}^{\frac{1}{2}}$ integrable, and also $S_{f^n}dz^2$ is $L_{loc}^{\frac{1}{2}}$ integrable.

%%%%%%
%%%%%%
\section{Local Convergence of $S_{f^n}$}\label{Local convergence}

%%%%%%
%%%%%%

%%%%%%
%%%%%%
In this section, our main goal is to prove  Theorem \ref{maintheorem1}, the local convergence of the normalized $S_{f^n}$.

\subsection{Bounded Fatou components}  \label{bounded Fatou component}

In this subsection, we are trying to show that $\lim_{n\to\infty}\frac{S_{f^n}}{d^{2n}}=0$
 on the bounded Fatou components for any degree $d\geq 2$ polynomial $f$, which is not conformally conjugate to $z^d$.
\begin{theorem}\label{bounded fatou converge}
For any $z\notin \Precrit(f)$ in a bounded Fatou component of a polynomial $f$ with degree $d\geq 2$ , which is not conformally conjugate to $z^d$, then we have:
$$\lim_{n\to\infty}\frac{S_{f^n}(z)}{d^{2n}}=0.$$
Moreover, this is a local uniform convergence.
\end{theorem}
\proof

First, assume that $z$ is attracted to some fix point $z_1$ (attracting or parabolic fix point), and $z_1$ is not a critical point. Then,
$$0<\lambda=|f'(z_1)|\leq 1$$
For any fixed $0<\epsilon<1$, since we have $f^n(z)$ converges to $z_1$, there exists $N_o\in\N$ and $M<\infty$, such that for any $n>N_o$, we have:
$$|f'(f^n(z))|\leq 1+\epsilon\ \textup{and}\ |S_f(f^n(z))|<M$$
By (\ref{schw equation}),
{\small
\begin{eqnarray*}
|S_{f^n}(z)|&=&|\sum_{i=1}^{n-1}S_f\circ f^{i}(z)\cdot((f^{i})'(z))^2+S_f(z)|\\
&\leq &\sum_{i=1}^{n-1}|S_f\circ f^{i}(z)\cdot((f^{i})'(z))^2|+|S_f(z)|
\end{eqnarray*}}
Since we have
$$|S_f(f^{n}(z))((f^{n})'(z))^2|<M\cdot M_1\cdot (1+\epsilon)^{2n},\ for\ any\ n>N_o.$$
where $M_1= |(f^{N_o})'(z)|$, by the fact that $1+\epsilon<d$ for $d\geq2$ and $\epsilon<1$, it is obvious that $\lim_{n\to\infty}\frac{S_{f^n}(z)}{d^{2n}}=0$
is satisfied.

\medskip

Second assume $z$ is attracted to a critical fix point $z_1$. Without loss of generality, we can assume $z_1=0$, so $f=az^r+bz^{r+1}+\cdots$, with $a\neq0$ and $2\leq r\leq d-1$. By Prop. \ref{basic prop}, we can study the Schwarzian limit at forward iterate of $z$ instead of the Swhwarzian limit at $z$. Then we can assume that $z$ is close to $0$. Let's conjugate $f$ to $z^r$ near $0$ by a conformal map $\phi$, such that $\phi(0)=0$ and $\phi'(0)\neq 0$,
$$\phi\circ f \circ \phi^{-1}=z^r$$
By the cocycle property of Schwarzian and Example 1,
\begin{eqnarray*}
S_{f^{n}}(z)&=&S_{\phi^{-1}\circ z^{r^n} \circ \phi}(z)\\
&=&S_{\phi^{-1}}\left((\phi(z))^{r^{n}}\right)\cdot \left((\phi^{-1})'(\phi(z)^{r^{n}})\cdot (r^n-1)\phi(z)^{r^{n}-1}\phi'(z)\right)^2\\
&\ &+\frac{1-r^{2n}}{2\phi(z)^2}\cdot (\phi'(z))^2 +S_\phi(z)
\end{eqnarray*}
Since $z$ is close to $0$ and $2\leq r \leq d-1$, then  $\lim_{n\to\infty}\frac{S_{f^n}(z)}{d^{2n}}=0$ is obviously true in this case by the above formula.

\medskip

Third, assume that $z$ is in some Siegel disc $\D_o$ fixed by $f$. Similar with the previous case, we can move the center of $\D_o$ to $0$, and conjugate $f$ by $\phi$ on this Siegel disc to a rotation map, i.e.
$$\phi\circ f =\lambda \cdot \phi, \textup{ with $|\lambda|=1$}$$
where we have $\phi(0)=0$ and $\phi'(0)=1$. And by (\ref{basic schw dicompositon}),
\begin{eqnarray*}
S_{f^{n}}(z)&=&S_{\phi^{-1}(\lambda^n \phi)}(z)\\
&=&S_{\phi^{-1}}(\lambda^n \phi(z))\cdot ((\phi^{-1})'(\lambda^n \phi(z)))^2 +S_\phi(z)
\end{eqnarray*}
Since $\lambda^n \phi(z)$ is in $\phi$'s image of some compact subset of $\D_o$ for any $n$, then $S_{\phi^{-1}}(\lambda^n \phi(z))$ is uniformly bounded.  So $\lim_{n\to\infty}\frac{S_{f^n}(z)}{d^{2n}}=0$ is obviously satisfied in this case by the above formula.

 \medskip
 From above arguments, it is not hard to see the convergence is local uniform. For points in the  periodic bounded Fatou components and not in $\Precrit(f)$, we can use similar arguments to show that the result of this theorem is satisfied. And for points attracted to  periodic bounded Fatou components, Prop. \ref{basic prop} shows that the result is satisfied too. \qed

 \bigskip

\subsection{Basin of infinity.} In this subsection, we are going to prove the local convergence of $S_{f^n}$ on the basin of infinity.

Consider the following escape-rate function of a polynomial $f$ with degree $d\geq 2$:
$$G_f(z) = \lim_{n\to\infty} \frac{1}{d^n} \log^+ |f^n(z)|,$$
where $\log^+|x|=\textup{max}(\log|x|,0)$. The escape-rate function $G_f(z)$ is the Green function of the basin of infinity. So it is harmonic on the basin of infinity. Actually, it is a subharmonic function on $\C$. By taking the partial derivative of $G_f(z)$, we get $g(z)=\frac{\del G_f(z)}{\del z}$ is a holomorphic function on
the basin of infinity; the zeros of $g(z)$ are exactly the points in $\Precrit(f)$.  Moreover, from the definition of $G_f(z)$, it is a limit of harmonic functions converging locally uniformly. The derivatives of this harmonic functions converge.  Then we know that the partial derivative commutes with the limit, i.e.

\begin{eqnarray}\label{g(z)}
g(z)&=&\frac{\del G_f(z)}{\del z}=\lim_{n\to\infty} \frac{\del\frac{1}{d^n} \log |f^n(z)|}{\del z}\\[6pt]\notag
&=&\lim_{n\to\infty} \frac{1}{2d^n} \frac{\del\log (f^n(z)\bar{f^n}(z))}{\del z}=\lim_{n\to\infty}\frac{1}{2} \frac{(f^n(z))'}{d^nf^n(z)}\notag
\end{eqnarray}

\begin{theorem}\label{basin of infinity converge}
For any $z_o\notin$ $\Precrit(f)$ in the basin of infinity of a polynomial $f$ with degree $d\geq 2$, we have:
$$\lim_{n\to\infty}\frac{S_{f^n}(z_o)}{d^{2n}}=-2(g(z_o))^2=-2\left(\frac{\del G_f(z_o)}{\del z}\right)^2$$
Moreover, this is a local uniform convergence.
\end{theorem}
Proof: Since $g(z)=0$ if and only if $z\in \Precrit(f)$, then we have $g(z_o)\neq 0$, because $z_o\notin\Precrit(f)$.

Let $r_n=\frac{(f^n(z_o))'}{2g(z_o)d^nf^n(z_o)}$, by (\ref{g(z)}), we get
$$\lim_{n\to\infty}r_n=1$$
Moreover, let $f(z)=a_{d}z^d+a_{d-1}z^{d-1}+\cdots+a_0$ with $a_{d}\neq 0$, and an easy calculation shows that
$$S_f(z)=\frac{1-d^2}{2z^2}h(z),\ and\ \lim_{z\to\infty} h(z)=1,$$
where $h(z)$ is a rational function.
Let $s_n=h(f^n(z_o))$, since $\lim_{n\to\infty} f^n(z_o)=\infty$, so we get $\lim_{n\to\infty} s_n=1$. Then
$$
\begin{array}{lll}
S_f(f^n(z_o))((f^n)'(z_o))^2
&=\frac{(1-d^2)s_n((f^n)'(z_o))^2}{2(f^n(z_o))^2}\\[6pt]
&=2(1-d^2)s_nr_n^2(g(z_o))^2d^{2n},

\end{array}
$$

Substituting the above formula into (\ref{schw equation}), we get
$$
\begin{array}{lll}
\frac{S_{f^n}(z_o)}{d^{2n}}
&=\frac{\sum_{i=1}^{n-1}S_f(f^{i}(z_o))((f^{i})'(z_o))^2+S_f(z_o)}{d^{2n}}\\[6pt]
&=2(g(z_o))^2(1-d^2)\frac{\sum_{i=0}^{n-1}s_ir_i^2d^{2i}}{d^{2n}}

\end{array}
$$

Because both $r_n$ and $s_n$ converge to $1$ as n tends to $\infty$, it follows that
$$
\begin{array}{lll}
\lim_{n\to \infty}\frac{S_{f^n}(z_o)}{d^{2n}}
&=2(g(z_o))^2(1-d^2) \lim_{n\to\infty} \frac{1-d^{2n}}{(1-d^2)d^{2n}} \\[6pt]
&=-2(g(z_o))^2=-2\left(\frac{\del G_f(z_o)}{\del z}\right)^2

\end{array}
$$

The fact that this convergence is local uniform can be deduced from the fact that both $r_n(z)$ and $s_n(z)$ converge locally uniformly.\qed

\medskip

{\bf Remark.}  Alternately, the result of Theorem \ref{basin of infinity converge} can be seen in the language of  ``Schwarzian between conformal metrics". On the complex plane $\C$, we define:
$$\widehat{S}(e^{\sigma_1}|dz|^2, e^{\sigma_2}|dz|^2)=\left(\sigma_{1zz}-\frac{1}{2}\sigma_{1z}^2-(\sigma_{2zz}-\frac{1}{2}\sigma_{2z}^2)\right)dz^2,$$
where $\sigma_{z}=\frac{\partial \sigma}{\partial z}$. Easy to know,  we have
 \begin{itemize}
\item $S_fdz^2=\widehat{S}(f^*|dz|^2, |dz|^2)$, for any non constant holomorphic map $f$.
\item $\widehat{S}(c_1\rho_1|dz|^2, c_2\rho_2|dz|^2)=\widehat{S}(\rho_1|dz|^2, \rho_2|dz|^2)$, for any positive constant numbers $c_1$ and $c_2$.
\item $\widehat{S}(\rho_1|dz|^2, \rho_3|dz|^2)=\widehat{S}(\rho_1|dz|^2, \rho_2|dz|^2)+\widehat{S}(\rho_2|dz|^2, \rho_3|dz|^2).$

\end{itemize}
Then we have

\includegraphics[width=145mm]{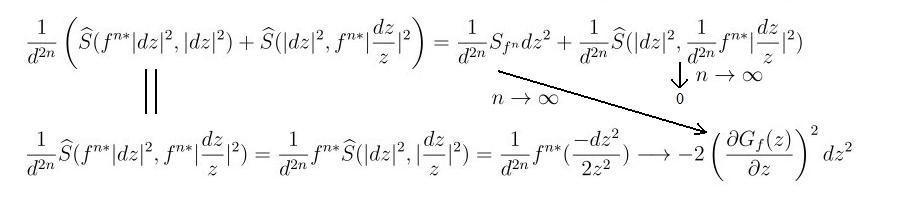}

For more about the Schwarzian of conformally equivalent Riemannian metrics, please refer to \cite{OS}.

\bigskip
{\bf Nonlinearity.} Similar with the Scharzian derivative, we can define nonlinearity $N_f$ of a nonconstant holomorphic function $f$ on the complex plane.
 $$N_f=\frac{f''}{f'}$$
Nonlinearity $N_f\equiv0$ if and only if $f$ is an affine transformation. Sometimes we can view $N_f$ as a one form $N_fdz$. We have the following cocycle property:
$$N_{f\circ g}dz=g^*(N_fdz)+N_gdz$$
Moreover, $N_fdz$ has a pole of order one at critical point of $f$. Using the same argument in Theorem \ref{maintheorem1}, we have
\begin{theorem}\label{nonlinearity}
Let $f$ be an polynomial with degree $d\geq 2$ and not conformally conjugate to $z^d$, and $X_o$ be its basin of infinity. Then we have:
\begin{itemize}
\item $\lim_{n\to \infty} \frac{(f^n)'dz}{d^nf^n}=\partial G_f$, on $X_o$.
\item $\lim_{n\to \infty} \frac{N_{f^n}dz}{d^n}=\partial G_f$, on $\C\setminus \overline{\Precrit (f)}$.
\item $\lim_{n\to \infty} \frac{(f^{n})'''dz^2}{d^{2n}(f^n)'}=-\frac{(\partial G_f)^2}{2}$, on $\C\setminus \overline{\Precrit (f)}$.
\end{itemize}
In each case, the convergence is locally uniform.
\end{theorem}

{\bf Proof of Theorem \ref{maintheorem1}.}
Since we have $G_f\equiv 0$ on the bounded Fatou components of $f$, then this theorem is an easy consequence of Theorem \ref{bounded fatou converge} and Theorem \ref{basin of infinity converge}.\qed

\bigskip

\subsection{Global convergence on the Fatou set.} Recall the definition of $L^{\frac{1}{2}}_{loc}$ integrability and $L^{\frac{1}{2}}_{loc}$ convergence of meromorphic quadratic differentials on the Riemann surface in \S  \ref{definition of quadratic differential}. For any rational function $f$ with degree $d\geq 2$, by Lemma \ref{L-1/2 integrable},  $S_{f^n}dz^2$ is $L^{\frac{1}{2}}_{loc}$ integrable.
\begin{cor}\label{maincor}
Let $S_{f^n}dz^2$ be the meromorphic quadratic differential on $\P^1$ determined by the Schwarzian derivative of $f^n$, where $f$ is a polynomial with degree $d\geq 2$ and not conformally conjugate to $z^d$. Then
$$\lim_{n\to\infty}\frac{S_{f^n}(z)}{d^{2n}}dz^2=-2\left(\frac{\del G_f(z)}{\del z}\right)^2dz^2,$$
on the Fatou set (including $\infty$) of $f$, converging in the sense of $L^{\frac{1}{2}}_{loc}$.

\end{cor}

\proof The proof of this corollary follows easily from the arguments in the proof of Theorem \ref{maintheorem1}, together with the triangle inequality. So we omit the details here.\qed

 %%%%%%
%%%%%%
\bigskip

\section{Geometric Limit of Metric Spaces} \label{Ultralimit of Metric Spaces}

In this section, we discuss the possible limit of metric spaces coming from Schwarzian derivatives.

\subsection{Ultrafilter and ultralimit}
%In this subsection, we describe the ultrafilter and then get an ultralimit of metric spaces coming from the Schwarzian derivatives of polynomial iterates.
 A non-principal ultrafilter $\omega$ is a set consisting of a collection of subsets of $\N$, satisfying
\begin{itemize}
\item If $A\subset B\subset \N$ and $A\in \omega$, then $B\in \omega$.
\item For any disjoint union $\N=A_1\cup A_2\cdots \cup A_n$, there exists one and only one $A_i\in \omega$.
\item  For any finite set $A\subset \N$, $A \notin \omega$,
\end{itemize}
\medskip
We can view a non-principal ultrafilter $\omega$ as a finitely additive measure on $\N$, only taking values in $\{0,1\}$, where any finite subset of $\N$ has measure zero and $\N$ has measure $1$.  By Zorn's lemma there exists some non-principal ultrafilter, and non-principal ultrafilter on $\N$ is not unique.  Hereafter, we fix a non-principal ultrafilter $\omega$. For more details about the ultrafilter, please refer to Chapter I \S 5 \cite{MA}.

\bigskip
Let $Y$ be a compact Hausdorff space. For any sequence of points $\{y_i\}_{i=1}^\infty\subset Y$, there is an unique point $y_o\in Y$ such that $\{i| y_i \in U\}\in \omega$ for any open set $U$ containing $y_o$. This $y_o=\lim_\omega y_i$ is denoted as the ultralimit  of $\{y_i\}_{i=1}^\infty$.

Let $\{(Y_n, d_n, b_n)\}_{n=1}^{\infty}$ be a sequence of pointed metric spaces. Let $\widetilde{Y}_\omega$ be the set of all the sequences $(y_n)$ with $y_n\in Y_n$ satisfying:
$$\lim_\omega d_n(y_n,b_n)< \infty;$$
here, the ultralimit is taken in the space of $[0,+\infty]$. Set
$$\widetilde{d}_\omega((x_n),(y_n)):=\lim_\omega d_n(x_n,y_n)<\infty,$$
with $(x_n)$ and $(y_n)\in \widetilde{Y}_\omega$. This is a pseudo-distance on $\widetilde{Y}_\omega$.
Let $(Y_\omega, d_\omega, b_\omega):=(\widetilde{Y}_\omega,\widetilde{d}_\omega,(b_n))/\sim$, where we identify points with zero $\widetilde{d}_\omega$-distance. The point metric space $(Y_\omega, d_\omega, b_\omega)$ is called the ultralimit of $\{(Y_n, d_n, b_n)\}_{n=1}^{\infty}$; see  Chapter I \S 5 \cite{MA}.

\medskip
\subsection{The ultralimit of the Schwarzian metrics}
Let $f$ be a polynomial with degree $d\geq 2$. The Schwarzian derivative $S_{f^n}$ determines a metric space $(X_n,d_n)$, where $X_n=\C\backslash \Crit(f^n)$ and $d_n$ is the arc length metric $ds=\sqrt{\left| \frac{S_{f^n}}{d^{2n}} \right| }|dz|$. This is a complete geodesic space with non positive curvature. Fix a point $a\notin \Precrit(f)$ on the basin of infinity. In this section we are considering the ultralimit $(X_{\omega},d_\omega,a_\omega)$ of the pointed metric spaces $(X_n, d_n, a)$.

\begin{prop}\label{prop of Xw space}
$(X_\omega,d_\omega,a_\omega)$ is a complete geodesic space.
\end{prop}
{\bf Proof:} The metric space $(X_\omega, d_\omega, a_\omega)$ is complete, since the ultralimit of any sequence of pointed metric spaces is complete; see   $\S1$ Lemma 5.53 \cite{MA}.  Moreover, the geodesic property is also preserved by passing to the ultralimit. The pointed metric space $(X_\omega, d_\omega, a_\omega)$ is  a geodesic space. Indeed, $\{(X_n, d_n, a)\}$ are pointed geodesic spaces. For any two points $x_\omega=(x_i)$ and $y_\omega=(y_i)$ in $X_\omega$, there is $(z_i)$ such that $d_i(x_i, z_i)=d_i(z_i, y_i)=\frac{1}{2}d_i(x_i,y_i)$. From this fact, we know that $z_\omega=(z_i) \in X_\omega$ and also
$$\lim_\omega \frac{1}{2}d_i(x_i, y_i)=\lim_\omega d_i(x_i, z_i)=\lim_\omega d_i(z_i, y_i),$$
i.e., $d_\omega(x_\omega, z_\omega)=d_\omega(z_\omega, y_\omega)=\frac{1}{2}d_\omega(x_\omega, y_\omega)$. So $x_\omega$ and $y_\omega$ have a midpoint in $X_\omega$. Which means $X_\omega$ is a geodesic space, since it is complete.\qed

\bigskip
\subsection{The flat metric on the basin of infinity}
%In this subsection, we are going to prove Theorem \ref{ultralimit spaces}. Lemma \ref{unique metric space} will be used in proving it.
Let $X_o$ be the basin of infinity of $f$ and $d_o$ be the arc length metric from $ds=\sqrt{2}\left|\frac{\partial G_f}{\partial z}\right||dz|$; see Figure \ref{flat metric}. The metric space $(X_o, d_o)$ is not complete. If the filled Julia set $K(f)$ of $f$ is disconnected, then we can complete $X_o$ as follows

\begin{lemma}\label{completion}
For polynomial $f$ with degree $d\geq 2$ and $K(f)$ disconnected, the metric completion $(\overline{X}_o, d_o)$ of $(X_o, d_o)$ is a quotient of $\C$, obtained by collapsing each connected component of $K(f)$ to a point. The completion $\overline{X}_o$ is  homeomorphic to $\R^2$. Moreover, $(\overline{X}_o, d_o)$ is a geodesic space.
\end{lemma}
\proof Since $K(f)$ is disconnected, the first two conclusions of this lemma follow immediately from the flat structure of the metric $d_o$; compare \cite{DM:trees}. A complete metric space is geodesic if there exists a midpoint for any two points on this metric space. It is obvious that $(\overline{X}_o,d_o )$ has this property.  So it is a geodesic space.\qed

\medskip

Let $E$ be the set of $\overline{X}_o \setminus X_o$. The set $E$ is totally disconnected. So we call $E$ the \textit{ends} of $X_o$.  Each point $e\in E$ corresponds to a connected component of $K(f)$. Let $C=X_o\cap \Precrit(f)$. For any \textit{end} in $E$, there is a sequence of annuli on $X_o\backslash C$ such that they nest down to this \textit{end} with $d_o$-diameters tending to zero; see the flat metric structure of $X_o$ in \cite{DM:trees}. The following lemma follows easily from Theorem \ref{maintheorem1}.

\begin{lemma}\label{curve converging}
For any piece wise smooth and compact curve $\gamma\subset X_o\setminus C$, the $d_n$-length $d_n(\gamma)$ converge to the $d_o$-length $d_o(\gamma)$. Moreover, for any point $p\in X_o\setminus C$, there is a small neighborhood $U\subset X_o\setminus C$ of $p$, such that $\lim_{n\to \infty} d_n(x,y)=d_o(x,y)$, uniformly for $x, y \in U$.
\end{lemma}

\bigskip
The following proposition is an essential ingredient in the proof of Theorem \ref{maintheorem1}.
\begin{prop}\label{unique metric space}
Let $f$ be a polynomial with degree $d\geq 2$ and $J(f)$ disconnected. For any geodesic metric $\widetilde{d}$ on $\overline{X}_o$ with $(X_o,\widetilde{d})$ locally isometric to $(X_o,d_o)$ under the identity map, we have $\widetilde{d}=d_o$ on $\overline{X}_o$.
\end{prop}

To prove this lemma, we need the tree structure of $(X_o,d_o)$; see \cite{DM:trees}. Specifically, this is a quotient $\pi: X_o \longrightarrow T(f)$, defined by collapsing each connected component of the level set of $G_f(z)$ to a point. There is a canonical map $F$ on $T(f)$ induced from $f$.
$$\xymatrix{ X_o \ar[d]_{\pi} \ar[r]^{f} &X_o  \ar[d]_{\pi}\\
		T(f)   \ar[r]^{F} &T(f)  }$$

The space $T(f)$ has a simplicial structure, and is equipped with a metric from $G_f$. The set of vertices is $V=\cup_{n,~m \in \Z}F^m(F^n(\textup{branch points}))$. The distance between two points in the same edge is given by the difference of their $G_f$ values. Let $S=\pi^{-1}(V)$. Then $X_o\setminus S$ consists of countably many connected components. Each of them is an annulus. So we view $X_o\setminus S$ as a set of annuli. The map $\pi$ is a one to one map from the annuli to the edges of $T(f)$. For any annulus $A\in X_o\setminus S$, it is a cylinder with finite height in the $d_o$-metric.

The height of $A$, denoted as $H(A)$, is equal to the length of the edge $\pi(A)$. Also, we define the level of $A$ as $L(A)$ to be the $G_f$ value of the middle point of $\pi(A)$. The closed geodesics of $A$ are the level sets of $G_f$. They have the same arc length in the $d_o$-metric, denoted as $C(A)$.

The map $f$ sends annulus in $X_o\setminus S$ to annulus. For each annulus $A\in X_o\setminus S$, there is a well defined local degree $d_A$, defined as the topological degree of $f|_A$. Moreover, let $N(f)=\textup{max}\{G_f(c)\ | c \textup{ is a critical point of f}\}$. We have the following properties:
\begin{itemize}
\item \begin{equation}\label{constant sum}
\sum_{B\in X_o\backslash S, L(B)=L(A)}C(B)=\sqrt{2}\pi
\end{equation}
\item \begin{equation}\label{degree map}
L(f(A))=d\cdot L(A), \ H(f(A))=d\cdot H(A),\  C(f(A))=\frac{d_A}{d}C(A)
\end{equation}
\item $d_A$ is equal to one plus the number of critical points (counted with multiplicity) enclosed in $A$.
\item The points enclosed in $A$ should have $G_f$ value less than $L(A)$. If $L(A)<N(f)$, then the annulus $A$ can not enclose the critical point(s) with $G_f$ value equals $N(f)$. And because $f$ has $d-1$ critical points counted with multiplicity, then $d_A\leq d-1$ when $L(A)<N(f)$. From (\ref{constant sum}) and (\ref{degree map}), for any annulus $A$ satisfying $d^nL(A)<N(f)$, we have $L(f^i(A))=d^iL(A)<N(f)$, i.e. $d_{f^i(A)}\leq d-1$ for any $1\leq i\leq n$. Consequently,
\begin{equation}\label{small c}
C(A)\leq (\frac{d-1}{d})^n\cdot C(f^n(A))\leq (\frac{d-1}{d})^n\cdot \sqrt{2}\pi
\end{equation}
\item Let $\left| f^{-1}(A) \right| $ be the number of annuli in the set $f^{-1}(A)$. We have $\left| f^{-1}(A) \right| $ equals to $d$ minus the total number of critical points (counted with multiplicity) enclosed in the annuli in $f^{-1}(A)$.
\item Let $A\in X_o\setminus S$. Any point $x\in \overline{X_o}$ with $G_f(x)\leq L(A)/d$ is enclosed in one of the annulus $B\in X_o\setminus S$ with $L(B)=L(A)$.
\item Any two annuli $A, B\in X_o\setminus S$ with $L(A)=L(B)$ have the same height.
\end{itemize}
 For more details about these annuli and the tree structure of the basin of infinity, please see \cite{DM:trees}.

\medskip
{\bf Proof of Proposition \ref{unique metric space}.} Since $X_o$ is dense in $\overline{X}_o$ for both the metrics $\widetilde{d}$ and $d_o$, then it suffices to prove that  $d_o(x,y)=\widetilde{d}(x,y)$ for any two points $x, y\in X_o$.
\medskip

One direction: $d_o(x,y)\geq \widetilde{d}(x,y)$. From the definition of $d_o$, for any $\epsilon >0$, there is an arc $\gamma \subset X_o$ connecting $x$ and $y$ with $d_o$ length $d_o(\gamma)<d_o(x,y)+\epsilon$. Since these two metrics are locally isometric on $X_o$, for the length $\widetilde{d}(\gamma)$ of $\gamma$ in the $\widetilde{d}$-metric, we have
$$\widetilde{d}(x,y) \leq \widetilde{d}(\gamma)=d_o(\gamma)<d_o(x,y)+\epsilon$$
Let $\epsilon\to 0$, then we get $\widetilde{d}(x,y)\leq d_o(x,y)$.

\medskip

The other direction: $d_o(x,y)\leq \widetilde{d}(x,y)$. Choose a geodesic $\widetilde{g}: [0,r] \mapsto \overline{X}_o$ in the $\widetilde{d}$-metric, with $\widetilde{g}(0)=x$ and $\widetilde{g}(r)=y$. Our goal is to construct an arc $\gamma \subset X_o$ from $\widetilde{g}$, connecting $x$ and $y$ with $d_o$ length $d_o(\gamma)<r+\epsilon$ for any fixed $\epsilon >0$.

From (\ref{small c}), there is some small level $l>0$ such that the set of annuli $T=\{ A \in X_o\backslash S\ | \ L(A)=l \}$ is not empty and $C(A)<\frac{\epsilon}{3}$ for any $A\in T$. Since all the annuli with the same level have the same height, so $h=H(A)$ for $A\in T$ is well defined. We order the elements of $T$ as $A_1, A_2,\cdots, A_k$, with $C(A_i)\leq C(A_j)$ for any $1\leq i<j\leq k$. Let $T_s=f^{-s}(T)$ and decompose $T_s$ into $T_s^-$ and $T_s^+$, with $T_s^+$ the set of annuli crossed by  $\widetilde{g}$. So the length of $\widetilde{g}$ in each annulus belonging to $T_s^+$ is at least the height of this annulus. Let $|\cdot|$ be the number of elements of a set. Then we have:
$$\sum_{s\in \N}|T^+_s|\cdot \frac{h}{d^s}\leq \textup{$\widetilde{d}$-length of $\widetilde{g}$}=r<\infty$$

Form the above formula, there is some big $n$, such that $|T^+_n|\cdot \frac{h}{d^n}<h$, i.e., $|T^+_n|<d^n$. As $|f^{-1}(A)|$ is equal to $d$ minus the total number of critical points (counted with multiplicity) enclosed in the annuli belonging to $f^{-1}(A)$, and $f$ has  $d-1$ critical points, then
$$|T_{s+1}|=|f^{-1}(T_s)|>d\cdot |T_s|-d, \ |f^{-n}(A)|\leq d^n$$
since any two annuli with the same level won't enclose each other, i.e., they won't share a common critical point inside. Consequently,
$$|T_{n}|>d^n\cdot |T|-d^n-d^{n-1}-\cdots-d>d^n\cdot k-2d^n$$
Combining the above formula with $|T^+_n|<d^n$, we have $|T^-_n|=|T_n|-|T^+_n|>d^n\cdot(k-3)$. From  (\ref{degree map}), for any $1\leq i\leq k$,
$$|f^{-n}(A_i)|\leq d^n, \ \frac{C(A_i)}{d^n}\leq C(B) \textup{ any } B\in f^{-n}(A_i)$$
We order $T_n$ as $B_1, B_2, \cdots, B_{|T_n|}$, such that $\beta_j=C(B_j)\leq \beta_{j+1}=C(B_{j+1})$ for any $1\leq j \leq |T_n|-1$. From the above formula, we have
 $$\beta_{d^n\cdot(i-1)+j}\geq \frac{C(A_i)}{d^n}, \textup{ $1\leq i\leq k-3$ and $1\leq j\leq d^n$,}$$
 Then we get,
%$$
%\begin{array}{324324}
%\sum_{B\in T^+_n}C(B)
%&=\sum_{B\in T_n} C(B)-\sum_{B\in T^-_n} C(B)\\[6pt]
%&\leq \sqrt{2}\pi-\sum_{i=1}_{k-3}\sum_{j=1}^{d^n}\beta_{d^n\cdot(i-1)+j}\\[6pt]
%&\leq \sum_{i=1}^{k}-\sum_{i=1}_{k-3}d^n\cdot \frac{C(A_i)}{d^n}\\[6pt]
%&=C(A_{k-2})+C(A_{k-1})+C(A_{k})\\[6pt]
%&\leq \epsilon.
%\end{array}
%$$
$$
\begin{array}{lll}
\sum_{B\in T^+_n}C(B)
&=\sum_{B\in T_n}C(B)-\sum_{B\in T^-_n}C(B)\\[6pt]
&\leq \sqrt{2}\pi-\sum_{i=1}^{k-3}\sum_{j=1}^{d^n}\beta_{d^n\cdot(i-1)+j}\\[6pt]
&\leq \sum_{i=1}^{k}C(A_i)-\sum_{i=1}^{k-3}d^n\cdot \frac{C(A_i)}{d^n}\\[6pt]
&=C(A_{k-2})+C(A_{k-1})+C(A_{k})\\[6pt]
&< \epsilon.
\end{array}
$$

Now, we can construct an arc $\gamma\subset X_o$ from $\widetilde{g}$ as follows (without loss of generality, we can always assume that both $G_f(x)$ and $G_f(y)$ are greater than $l$, and $n>1$):
\begin{itemize}
\item Choose a minimal $r_1\in [0,r]$ such that $\widetilde{g}(r_1)$ lying on the outer boundary of some annulus $B^1\in T^+_n$.
\item Choose a maximal $r_1'\in [0,r]$ such that $\widetilde{g}(r_1')$ lying on the outer boundary of the annulus $B^1$.
\item Replace the arc $\widetilde{g}([r_1, r_1'])\subset \widetilde{g}$ with a shortest curve $\gamma_1$ on the outer boundary of $B^1$ connecting $\widetilde{g}(r_1)$ and $\widetilde{g}(r_1')$.
\item Do the same thing as the previous three steps, we can find a minimal $r_2$ and maximal $r_2'$ in $[r_1',r]$, such that $\widetilde{g}(r_2)$ and $\widetilde{g}(r_2')$ lying on the outer boundary of some $B^2\in T^+_n$ and $r_2'-r_2>0$. Replace $\widetilde{g}([r_2, r_2'])\subset \widetilde{g}$ with a shortest curve $\gamma_2$ on the outer boundary of $B^2$ connecting $\widetilde{g}(r_2)$ and $\widetilde{g}(r_2')$.
\item Keep doing the same thing as previous step by step, we can replace sub-arcs of $\widetilde{g}$ by arcs $\gamma_i$ on the boundary of $B^i\in T^+_n$. This process will stop under finite steps, since $B^i\neq B^j$ for any $i<j$ and $|T^+_n|<\infty$.
\end{itemize}
From the above construction, we get a new arc $\gamma$ from $\widetilde{g}$. We have $\gamma\subset X_o$. In fact, for any point $p\in \gamma$, $p$ is not enclosed by any annulus $B\in T_n$. So we have $G_f(p)>\frac{h}{d^{n+1}}$. And
$$
\begin{array}{lll}
d_o(x,y)
&\leq d_o(\gamma)\leq r+d_o(\gamma_1)+d_o(\gamma_2)+\cdots\\[6pt]
&\leq r+C(B^1)+C(B^2)+\cdots\\[6pt]
&\leq r+ \sum_{B\in T^+_n}C(B)< r+\epsilon.
\end{array}
$$
Let $\epsilon\to 0$, we get $d_o(x,y)\leq r=\widetilde{d}(x,y)$.\qed

\bigskip

{\bf Proof of Theorem \ref{ultralimit spaces}.}
First, we construct a natural map $\rho$
  $$\rho : X_o \backslash C\mapsto X_\omega,$$
as: for any $p\in X_o\backslash C$, $\rho(p)=p_\omega=(p_n=p)\in X_\omega$. This map is well defined, since $\{d_n(p,a)\}$ is uniformly bounded by Lemma \ref{curve converging}. Let $p,\ q$ be two points in $X_o\backslash C$, and $p_\omega,\ q_\omega$ the $\rho$ image of $p,\ q$ in $X_\omega$. We want to show that $d_\omega(p_\omega,q_\omega)\leq d_o(p,q)$. For any $\epsilon>0$, we can choose a smooth curve on the basin of infinity with $d_o$-length less than $d_o(p,q)+\epsilon$. By a small perturbation, we can assume this curve does not pass any point in $C$, and this curve is also on the basin of infinity. Since this curve is compact, by Lemma \ref{curve converging}, the $d_n$-length of this curve converges to the $d_o$-length of this curve. So for $n$ big enough, the $d_n$-length of this curve is less than $d_o(p,q)+2\epsilon$. Then $d_\omega(p_\omega,q_\omega)\leq d_o(p,q)+2\epsilon$. Let $\epsilon\rightarrow 0$, we get $d_\omega(p_\omega,q_\omega)\leq d_o(p,q)$. Also, from Lemma \ref{curve converging}, this is a locally isometric and distance non-increasing embedding. Then we can extend the map $\rho$ from  $ X_o \backslash C$ to $ \overline{X_o}$:
$$\rho: \overline{X_o} \longrightarrow X_\omega$$

For any \textit{end} $e\in E$, let $K_e$ be the corresponding connected component of $K(f)$ and $X^e_\omega=\{(x_i)\in X_\omega|\lim_\omega x_i\in K_e\}$. And for any $c\in C$, let $X_\omega^c=\{(x_i)\in X_\omega|\lim_\omega x_i=c\}$. Obviously, the set
$$X_\omega=\left(\cup_{\alpha\in C\cup E}X^\alpha_\omega\right)\cup \rho(X_o\backslash C),$$
is a disjoint union, i.e.
$$X_\omega^{\alpha_1}\cap X_\omega^{\alpha_2}=\emptyset \textup{ and }X_\omega^{\alpha_1}\cap \rho(X_o\backslash C)=\emptyset\textup{, for any $\alpha_1 \neq \alpha_2 \in C\cup E$}$$
Indeed, there is a closed annulus $A$ in $X_o\backslash C$ with the corresponding parts of $\alpha_1$ and $\alpha_2$ in the two different components of $\C\setminus A$. Let $h$ be the distance between the two boundaries of $A$ in the $d_o$-metric. we have $h>0$. Since any arc connecting two points in distinct components of $\C\setminus A$ should across $A$. By Lemma \ref{curve converging}, the $d_n$-distance of the two boundaries of $A$ converge to $h$. So we have $d_\omega(X_\omega^{\alpha_1},X_\omega^{\alpha_2})\geq h>0$. Consequently, $X_\omega^{\alpha_1}\cap X_\omega^{\alpha_2}=\emptyset$. Similarly, any point $x_\omega \in \rho(X_o\backslash C)$, we have $d_\omega(X_\omega^{\alpha_1},x_\omega )>0$, then $X_\omega^{\alpha_1}\cap \rho(X_o\backslash C)=\emptyset$.

From above, for any two distinct points $\alpha$, $\beta$ in $C\cup E$ and $p_\omega\in \rho(X_o\backslash C)$, we have $d_\omega(\rho(\alpha),X_\omega^\beta)>0$ and $d_\omega(\rho(\alpha),p_\omega)>0$. Moveover, since $\overline{X}_o$ is homeomorphic to $\R^2$ by Lemma \ref{completion}, so $\rho :\overline{X}_o\mapsto \rho(\overline{X}_o)$ is a distance non-increasing homeomorphism. We want to show $(\rho(\overline{X}_o),d_\omega)$ is a geodesic space and $\rho$ is locally isometric at the points in $C$.  Then by Proposition
\ref{unique metric space}, we may conclude that $\rho$ is an isometric embedding.

For any $x_\omega \in X_\omega^\alpha$ and $y_\omega \notin X_\omega^\alpha$, with $\alpha \in C\cup E$, we want to show any geodesic $g_\omega$ connecting $x_\omega$ and $y_\omega$ should pass through $\rho(\alpha)$. For any closed annulus $A\subset X_o\backslash C$ with the points $\lim_\omega (x_i)$ and $\lim_\omega (y_i)$ lying in the different components of $\C\setminus A$, as previous, $X_\omega\setminus \rho(A)$ consists two connected components, with distance at least the distance of the two boundaries of $A$ in the $d_o$-metric. And since $x_\omega$ and $y_\omega$ are in different components of  $X_\omega\setminus \rho(A)$, the geodesic connecting them should intersect $\rho(A)$. Let $\{A_i\}_{i=1}^{\infty}$ be a sequence of such annuli, and they nest down to  $\alpha$ in the sense that $\lim_{i\to \infty}$diameter$(\alpha \cup A_i)=0$ in the $d_o$-metric.  Choose some $x_\omega^i$ in $\rho(A_i)\cap g_\omega$. Since $\{A_i\}_{i=1}^{\infty}$ nest down to  $\alpha$ and the map $\rho$ does not increase the distance, then we have that $\{x_\omega^i\}$ converges to $\rho(\alpha)$. As the geodesic is compact, we know that $\rho(\alpha)$ should in the geodesic $g_\omega$.

If there is a geodesic $g_\omega\subset X_\omega$ with two end points in $\rho(\overline{X}_o)$ such that it has some point $p_\omega \in g_\omega \cap (X_\omega \setminus \rho(\overline{X}_o))$. Assume $p_\omega$ belongs to $X^\alpha_\omega$ with $\alpha \in C\cup E$. Then $p_\omega$ divides $g_\omega$ in to two parts. Each of these two parts should pass though $\rho(\alpha)$. Which means $g_\omega$ can not be the shortest curve connecting the two end points. In all, any geodesic with two ends in $\rho(\overline{X}_o)$ should be contained in $\rho(\overline{X}_o)$. So $(\rho(\overline{X}_o), d_\omega)$ is a geodesic space, since $X_\omega$ is a geodesic space. Plus $\rho|_{X_o\backslash C}$ is locally isometric and distance non-increasing embedding, we get $\rho|_{X_o}$ is locally isometric.

\medskip
For any $\alpha \in E\cup C$, since there always exists a sequence of $\{c_i\}_{i=1}^\infty \subset C$ converging to $\alpha$, then for any $l>0$ big enough, we can always choose a sequence of points $\{x_i\}_{i=1}^{\infty} \subset X_o\backslash C$, such that $d_i(x_i, a_i)=l$ and $\lim_\omega x_i$ in $\alpha$'s corresponding subset of $\C$. Then we have $x_\omega=(x_i)\in X_\omega^\alpha$ and $d_\omega(x_\omega, a_\omega)=l$.

For any $c\in C$, to prove that $X^c_\omega$ is a real ray, it suffices to prove that for any two sequence $\{x_n\}$ and $\{y_n\}$ converging to $c$, with
$$\lim_{n\to \infty} \left[d_n(x_n,a)=d_n(y_n,a)\right]=l>l_o=d_o(a,c)$$
then, we have $\lim_{n\to \infty} d_n(x_n,y_n)=0$. This follows easily from Lemma \ref{basic lemma} proved below.\qed

\bigskip

\subsection{Real rays attached to $C$} In this subsection, we are going to complete the proof of Theorem \ref{ultralimit spaces} by showing Lemma \ref{basic lemma}. What remains is to show that the extra pieces of $X_\omega$ are real rays attached to $X_o$. The basic idea is to show that $X_n$ has no ``bulb" near the critical points of $f^n$. For doing this, we need to use the fact that $X_n$ is a metric space with non-positive curvature.

Let $S$ be a closed Riemann surface with genus $g\geq 2$ and $Q$ be a holomorphic quadratic differential on $S$. Then $Q$ determines a flat metric $ds^2=|Q|$ with finite singularities on $S$ at zeros of $Q$. At non-singular points, it is flat, so it has curvature $0$; at singular points, it's a cone with angle $k\pi$ for $3\leq k\in \N$. So at the singular points, it has negative curvature. Then in this metric, $S$ is a complete arc length metric space with non-positive curvature. For the definition and properties of the curvature, please refer to p. 159 \cite{MA}. Lift this metric to the universal cover $\widetilde{S}$ of $S$, we also get a metric on $\widetilde{S}$ with non-positive curvature.

\begin{lemma}\label{complete universalcover}
$\widetilde{S}$ is a complete CAT$(0)$ unique geodesic space, any geodesic locally is a straight line at non singular point.
\end{lemma}
\proof By \cite{LA}, $\widetilde{S}$  is an unique geodesic space (any two points are connected by an unique geodesic) with geodesic locally straight line, and $\widetilde{S}$ is also complete. Moreover, since it is a complete and simply connected metric space with non positive curvature, by Cartan-Hadamard Theorem, such space is a CAT$(0)$ space; see p. $193$ \cite{MA}.\qed

\bigskip

Fix a point $c\in C$ and some very small $\epsilon>0$. Let $c_\epsilon\subset X_o\backslash C$ be the closed curve at $\epsilon$-distance from $c$ in the $d_o$-metric. For each $n$, choose some closed geodesic (in the $d_n$-metric) $c_n\subset X_n$ such that $c_n$ is sufficiently close to $c$; see Lemma \ref{length of normalized trajectory}. By Lemma \ref{curve converging} and Lemma \ref{length of normalized trajectory}, there is some $M<\infty$, such that the $d_n$-length of $c_\epsilon$ $d_n(c_\epsilon)<M\cdot \epsilon$ and $d_n(c_n)<M/d^n$. Choose $x_n\in c_\epsilon$ and $y_n \in c_n$, such that $r_n=d_n(x_n, y_n)=d_n(c_\epsilon,c_n)$. We can do this is because both $c_\epsilon$ and $c_n$ are compact. Choose a geodesic $g_n:[0,r_n]\mapsto X_n$ with $g_n(0)=x_n$ and $g_n(r_n)=y_n$. Let $A_n$ be the open annulus bounded by $c_\epsilon$ and $c_n$. We have that $g_n((0,r_n))\subset A_n$. Otherwise it won't be the shortest curve connecting the two boundaries of $A_n$.

We can construct closed Riemann surface $S_n$ with genus $g\geq 2$ from $X_n$. The metric space ($X_n,d_n$) has finitely many infinite cylinders (cylinder with infinite height). Each such infinite cylinder lying in some neighborhood of a critical point (including $\{\infty \}$) of $f^n$. So cut the infinite cylinders off $X_n$ along some of the closed geodesics inside the cylinders, such that the closed geodesic is much closer than $c_n$ to the critical point. Then we get $X_n'$. Double $X_n'$, and glue them together along the corresponding boundaries to get a closed surface $S_n$. The metric on $S_n$ is the obvious metric induced from $X_n'$. Consider the universal cover $\widetilde{S}_n$ of $S_n$ with the induced metric $\widetilde{d}_n$ from $S_n$. Topologically, $\widetilde{S}_n$ is a unit disk. Let $\widetilde{A}_n\subset \widetilde{S}_n$ be one of the connected components of the preimage of $A_n\subset X_n'$. Then $\widetilde{A}_n$ is a strip on $\widetilde{S}_n$ separating $\widetilde{S}_n$ into two connected components. Since the projection of any curve connecting this two components should be some curve across $A_n\subset X_n'$, the distance of these two components is $d_n(x_n,y_n)$ obtained by some lift $\widetilde{g}_n$ of $g_n$ connecting these two components.

The boundaries $\partial\widetilde{A}_n$ of $\widetilde{A}_n$ are two curves in the preimages of $c_\epsilon,c_n\subset X_n'$. The preimage of $g_n$ on $\widetilde{A}_n$ cuts $\widetilde{A}_n$ into quadrilaterals. All of them can be mapped into each other by some isometry of $\widetilde{S}_n$. Pick one of these quadrilateral $\widetilde{B}_n$ with $\widetilde{g}_n\subset \partial\widetilde{B}_n$. Then $\widetilde{B}_n$ is a copy of the lift of $A_n\backslash g_n$. Denote $\widetilde{c}_\epsilon$ and  $\widetilde{c}_n$ as the lift of $c_\epsilon$ and $c_n$ on $\partial\widetilde{B}_n$, and $\widetilde{g}_n^o$ the other lift of $g_n$ on $\partial\widetilde{B}_n$.

\begin{lemma}\label{non bulgs}
For any two points $\widetilde{z}_1\in \widetilde{c}_\epsilon $ and $\widetilde{z}_2\in \widetilde{c}_n$, there is an unique  geodesic $\widetilde{g}$ connecting these two points, and if we varies  $\widetilde{z}_1$ and $\widetilde{z}_2$ continuously, then $\widetilde{g}$ varies continuously. In particular, any point $\widetilde{q}\in \widetilde{B}_n$, there is some geodesic $\widetilde{g}_q$ passing though $\widetilde{q}$ with two ends in $\widetilde{c}_\epsilon$ and $\widetilde{c}_n$.
\end{lemma}
\proof By Lemma \ref{complete universalcover}, there is an unique geodesic $\widetilde{g}$ connecting $\widetilde{z}_1$ and $\widetilde{z}_2$. Because $\widetilde{S}_n$ is a complete and simply connected CAT$(0)$ metric space, by  Cartan-Hadamard theorem in p. 193 \cite{MA}, geodesic varies continually with respect to the two end points.

Assume there is some point $\widetilde{q}\in \widetilde{B}_n$ such that any geodesic with two ends in $\widetilde{c}_\epsilon $ and $ \widetilde{c}_n$ won't pass though it. Choose $\widetilde{z}_1(t)\in \widetilde{c}_\epsilon $ and $\widetilde{z}_2(t)\in \widetilde{c}_n$ varies from the ends of $\widetilde{g}_n$ to the ends of $\widetilde{g}_n^o$. Then the corresponding geodesics varies from $\widetilde{g}_n$ to $\widetilde{g}_n^o$ without touching $\widetilde{q}$. From this we get that, in $\widetilde{S}_n\backslash \widetilde{q}$, $\partial\widetilde{B}_n$ is homotopic to a point. This is impossible since  $\widetilde{S}_n$ is topologically a disc.\qed

%\begin{figure}
%\includegraphics[width=6in]{mypdf.pdf}
%\caption{Here is a picture.}
%\label{my figure}
%\end{figure}

\begin{lemma}\label{basic lemma}
Let $\widetilde{q}_1$ and $\widetilde{q}_2$ be two points in $\overline{\widetilde{B}}_n\backslash(\widetilde{c}_n\cup\widetilde{c}_\epsilon)$, with $l_1=\widetilde{d}_n(\widetilde{q}_1, \widetilde{c}_\epsilon)$ and $l_2=\widetilde{d}_n(\widetilde{q}_2, \widetilde{c}_\epsilon)$. Then $\widetilde{d}_n(\widetilde{q}_1,\widetilde{q}_2)\leq |l_1-l_2|+9r_3+7r_4$, where $r_3$ and $r_4$ are $\widetilde{d}_n$-lengths of $\widetilde{c}_\epsilon$ and $\widetilde{c}_n$.

\end{lemma}

\proof As in Lemma \ref{non bulgs}, we can choose geodesics $\widetilde{g}_i: [0,r_i]\to \widetilde{S}_n$ of length $r_i$ passing though $\widetilde{q}_i$, with $\widetilde{g}_i(0)\in \widetilde{c}_\epsilon$ and  $\widetilde{g}_i(r_i)\in \widetilde{c}_n$ for $i=1,2$. Also, we have geodesic $\widetilde{g}_o:[0,r_o]\to \widetilde{S}_n$ with $\widetilde{g}_o(0)=\widetilde{g}_1(0)$ and $\widetilde{g}_o(r_o)=\widetilde{g}_2(r_2)$.

 %$s_1, s_2\in [0,1]$ such that the geodesic $\widetilde{g}_{s_i}$ passing $\widetilde{q}_i$. Denote $\widetilde{g}_{s_o}$ as the geodesic connecting $\widetilde{z}_1(s_1)$ and $\widetilde{z}_2(s_2)$.

For $1\leq i \leq 2$, there is $r_i'$ such that $\widetilde{q}_i=\widetilde{g}_{s_i}(r_i')$ with $0<r_i'<r_i$. And since $l_i=\widetilde{d}_n(\widetilde{p}_i, \widetilde{c}_\epsilon)$ and $r_i'=\widetilde{d}_n(\widetilde{g}_i(0),\widetilde{p}_i)$, then $l_i\leq r_i'\leq l_i+r_3$.

First, assume that we have $r_o\leq r_1\leq r_2$. In the isosceles triangle with three vertices $\widetilde{g}_1(0), \widetilde{g}_2(r_2)$ and $\widetilde{g}_{1}(r_o)$, since $r_1'=\widetilde{d}_n(\widetilde{q}_1,\widetilde{g}_1(0))=\widetilde{d}_n(\widetilde{g}_{o}(r_1'),\widetilde{g}_1(0))$, by CAT$(0)$ property of $\widetilde{S}_n$, we have
$$\widetilde{d}_n(\widetilde{q}_1,\widetilde{g}_{o}(r_1'))\leq \widetilde{d}_n(\widetilde{g}_{1}(r_o),\widetilde{g}_2(r_2))$$
$$\leq \widetilde{d}_n(\widetilde{g}_{1}(r_o),\widetilde{g}_1(r_1))+\widetilde{d}_n(\widetilde{g}_1(r_1),\widetilde{g}_2(r_2))\leq (r_1-r_o)+r_4$$
In the isosceles triangle with three vertices $\widetilde{g}_2(r_2), \widetilde{g}_1(0)$ and $\widetilde{g}_{2}(r_2-r_o)$, since $r_2'-(r_2-r_o)=\widetilde{d}_n(\widetilde{q}_2,\widetilde{g}_{2}(r_2-r_o))=\widetilde{d}_n(\widetilde{g}_{o}(r_2'-(r_2-r_o)),\widetilde{g}_1(0))$, by CAT$(0)$ property of $\widetilde{S}_n$, we have
$$\widetilde{d}_n(\widetilde{q}_2,\widetilde{g}_{o}(r_2'-(r_2-r_o)))\leq \widetilde{d}_n(\widetilde{g}_1(0),\widetilde{g}_{2}(r_2-r_o))$$
$$\leq \widetilde{d}_n(\widetilde{g}_1(0),\widetilde{g}_2(0))+\widetilde{d}_n(\widetilde{g}_2(0),\widetilde{g}_{2}(r_2-r_o))\leq r_3+(r_2-r_o)$$
Moreover, since $r_3$ and $r_4$ are $\widetilde{d}_n$-lengths of $\widetilde{c}_\epsilon$ and $\widetilde{c}_n$, so we have
$$ |r_1-r_o|\leq r_3+r_4 \textup{ and } |r_2-r_o|\leq r_3+r_4$$
Consequently,
$$
\begin{array}{lll}
\widetilde{d}_n(\widetilde{q}_1,\widetilde{q}_2)
&\leq\widetilde{d}_n(\widetilde{q}_1,\widetilde{g}_{o}(r_1'))+\widetilde{d}_n(\widetilde{g}_{o}(r_1'),\widetilde{g}_{o}(r_2'-(r_2-r_o)))+\widetilde{d}_n(\widetilde{g}_{o}(r_2'-(r_2-r_o)),\widetilde{q}_2)\\[6pt]
&\leq((r_1-r_o)+r_4)+|(r_2'-(r_2-r_o))-r_1'|+(r_3+(r_2-r_o))\\[6pt]
&\leq|r_1-r_o|+|r_2'-r_1'|+2|r_2-r_o|+r_3+r_4\\[6pt]
&\leq 2(r_3+r_4)+(2r_3+|l_1-l_2|)+2\cdot2(r_3+r_4)+r_3+r_4\\[6pt]
&=|l_1-l_2|+9r_3+7r_4,
\end{array}
$$

Second, for all other cases, similarly, we can always get:
$$\widetilde{d}_n(\widetilde{q}_1,\widetilde{q}_2)\leq |l_1-l_2|+9r_3+7r_4.$$\qed

\bigskip\bigskip
\def\cprime{$'$}

%\bibliographystyle{../tex/bib/math}
%\bibliography{../tex/bib/math}
%\bibliographystyle{math}
%\bibliography{refs}

 \end{document}